\documentclass[preprint,authoryear,12pt]{elsarticle}

\renewcommand{\baselinestretch}{1.25}
\usepackage{mathtools}
\usepackage[normalem]{ulem}
\usepackage{amsmath}
\usepackage{amssymb}
\usepackage{amsthm}
\usepackage[T1]{fontenc}
\usepackage[latin1]{}
\usepackage[utf8]{inputenc}
\usepackage{amsmath,amsfonts,amscd,bezier}
\usepackage{threeparttable}
\usepackage{tabulary}
\usepackage{booktabs}				
\usepackage{multirow}
\usepackage{algpseudocode,algorithm}
\usepackage{enumitem}
\usepackage{caption}
\usepackage{subcaption}\usepackage{lmodern}
\usepackage{mathrsfs}
\usepackage{graphics, graphicx,color}
\usepackage{epsfig}
\usepackage{lscape}
\usepackage{natbib}

\usepackage{adjustbox}

\usepackage{float} 
\usepackage{rotating}

\newcounter{example}[section]

\newcounter{remark}[section]

\newcommand{\n}{\noindent}

\newcommand{\Expec}{\mathbb{E}}

\usepackage{verbatim} 
\usepackage{comment}  
\usepackage{hyperref}
\hypersetup{
		colorlinks=true, % false: boxed links; true:  colored links
    	linkcolor=blue,          	% color of internal links
    	citecolor=blue,        		% color of links to bibliography
    	filecolor=magenta,      		% color of file links
		urlcolor=blue,
		bookmarksdepth=4}

\journal{arXiv}

\begin{document}
\setlength{\abovedisplayskip}{0pt}
\setlength{\belowdisplayskip}{0pt}
\setlength{\abovedisplayshortskip}{0pt}
\setlength{\belowdisplayshortskip}{0pt}

\begin{frontmatter}

\title{Revisiting Gini for equitable humanitarian logistics}

\vspace{0cm}

\author[label1]{Douglas Alem\corref{cor1}}\ead{Douglas.Alem@ed.ac.uk (DA)}
\author[label1]{Aakil M. Caunhye}\ead{Aakil.Caunhye@ed.ac.uk (AMC)}
\author[label2]{Alfredo Moreno}\ead{alfredo-daniel.moreno-arteaga@hec.ca}

\address[label1]{University of Edinburgh Business School\\29 Buccleuch Place, EH8 9JS, Edinburgh, United Kingdom}

\address[label2]{Department of Logistics and Operations Management, HEC, Montreal\\ 3000 Chemin de la Cote-Sainte-Catherine, Canada}

\cortext[cor1]{Corresponding author}

\renewcommand{\baselinestretch}{1.25}

\newpage

\begin{abstract}
Modeling equity in the allocation of scarce resources is a fast-growing concern in the humanitarian logistics field. The Gini coefficient is one of the most widely recognized measures of inequity and it was originally characterized by means of the Lorenz curve, which is a mathematical function that links the cumulative share of income to rank-ordered groups in a population. So far, humanitarian logistics models that have approached equity using the Gini coefficient do not actually optimize its original formulation, but use alternative definitions that do not necessarily replicate that original Gini measure. In this paper, we derive the original Gini coefficient via the Lorenz curve to optimize the effectiveness-equity trade-off in a humanitarian location-allocation problem. We also propose new valid inequalities based on an upper-bounding Lorenz curve to tighten the linear relaxation of our model and develop a clustering-based construction of the Lorenz curve that requires fewer additional constraints and variables than the original one. The computational study, based on the floods and landslides in Rio de Janeiro state, Brazil, reveals that while alternative Gini definitions have interesting properties, they can generate vastly different decisions compared to the original Gini coefficient. In addition, viewed from the perspective of the original Gini coefficient, these decisions can be significantly less equitable.
\end{abstract}

\vspace{0.5cm}

%\begin{keyword}
%Gini; equity/fairness; humanitarian logistics; location-allocation; disaster relief; Lorenz curve; clustering.
%\end{keyword}

\end{frontmatter}

\newpage

\section{Introduction}
\label{sec:Introduction}

Even though equity in allocation of scarce public resources has been a concern of academics and practitioners for decades, there is no consensus on how to mathematically model this concept within humanitarian logistics optimization. Ultimately, the concept of equity \textit{per se} is subject to different interpretations and is very much context-dependent \citep{donmez2021humanitarian}. However, scholars have identified three key aspects that can guide decision makers and reduce their subjective biases in choosing an appropriate equity modelling approach: identifying an equity metric, constructing an equity function, and formulating an equity objective. An equity metric refers to the quantity, such as the coverage of victim needs, that the decision maker seeks to balance; an equity function calculates the equity level achieved via a given set of decisions (such as the amount of resources allocated to cover victim needs); and an equity objective relates to a given combination of metric and function that the decision maker seeks to maximize/minimize \citep{matl2018workload}. One of the oldest and most popular equity (fairness) approaches is the so-called Rawlsian approach \citep{rawls1973theory}, which in optimization problems, translates into min-max equity functions for optimizing the \textit{worst-off} entity or outcome (see \citet{manopiniwes2017stochastic, noyan2018stochastic,Aslan2019,rodriguez2020shortage}, and references therein).

More recently, there has been an increased interest in the popular Gini coefficient as an equity function for humanitarian operations. Already one of the ``most widely recognized (measure of inequity) in the economic and social welfare literature'' \citep{marsh1994equity, leclerc2012modeling}, the Gini coefficient has been recently employed in the study of a variety of humanitarian logistics problems such as pickup and distribution \citep{eisenhandler2019}, shelter location \citep{Mostajabdaveh2019}, location and distribution \citep{park2020supply}, vaccine allocation \citep{enayati2020optimal}, and joint prepositioning of emergency supplies \citep{rodriguez2020cost}. Modelling equity by means of the Gini coefficient is more easily interpretable by policymakers and obeys desirable theoretical properties, such as the Pigou–Dalton Principle of Transfers, which states that a positive transfer of income from a `richer' to a `poorer' individual should decrease the extent of income inequality and lead to a more equitable outcome if that transfer preserves the rank-order of incomes \citep{mehran1976linear,karsu2015inequity}.

The Gini coefficient was originally characterized by means of the Lorenz curve, which is a mathematical function that links the cumulative share of income to rank-ordered groups in a population. With the Lorenz curve characterization, the coefficient not only preserves the Pigou-Dalton Principle of Transfers, but also becomes a scale-invariant and relative measure, taking on values from 0 (perfect equity) and  1 (perfect inequity). While the Gini coefficient is readily computable on data, its use in decision-making models has been superseded by easier-to-implement alternative definitions of Gini, mainly relying on mean difference measures. One popular example is $\sum^{n}_{i,j}|x_i-x_j|/2n^2\bar{x}$, where $x_1,x_2,\dots,x_n$ are possible outcomes, and $\bar{x}$ is the average value \citep{Mostajabdaveh2019}. Another example is the measure proposed by \cite{mandell1991}, which is $\sum_{i,j>i}|q_jS_i-q_iS_j|/\sum_iS_i$, where $q_i$ is the proportion of individuals in group $i$ and $S_i$ is the quantity of resources allocated to group $i$. Although it is decidedly true that these measures maintain most of the interesting properties of equity measures, such as the Principle of Transfers, it is arguable how closely they are able to match the original formulation of the Gini coefficient via the Lorenz curve. 

Interestingly, all the aforementioned humanitarian logistics models that have approached equity using the Gini coefficient do not actually optimize the original formulation based on the Lorenz curve, but use alternative definitions that can be more straightforwardly embedded in optimization models. In this paper, we derive the original formulation of the Gini coefficient and embed it in a new objective function focused on the \textit{desideratum} effectiveness-equity trade-off for a class of humanitarian location-allocation problem. The effectiveness function determines the extent to which victim needs are covered, whereas the equity function measures the extent to which relief aid is fairly allocated amongst disaster-prone areas.

The contributions of this study are fourfold: 1) we derive the original Gini coefficient via the Lorenz curve in a decision-making setting using a mixed-integer combination of constraints and variables that can be easily embedded in any optimization model, 2) we propose new valid inequalities based on an upper-bounding, suboptimal Lorenz curve to tighten the linear relaxation of our model under the original Gini coefficient formulation and thus improve the numerical efficiency of our solution, 3) we propose a clustering-based construction of the Lorenz curve that requires fewer additional constraints and variables than the original Lorenz curve, and 4) we apply our models to a realistic location-allocation problem involving floods and landslides in Rio de Janeiro state, Brazil, and show that while alternative Gini definitions based on mean difference measures have interesting properties, they do not replicate the original Gini coefficient based on the Lorenz curve. More importantly, in a decision-making setting, the mean difference can generate vastly different decisions compared to the original Gini coefficient and in addition, viewed from the perspective of the original Gini coefficients, these decisions can be significantly less equitable. We would like to emphasize that although we develop the original Gini coefficient for a humanitarian logistics application, our formulation can be easily translated to any other applications, without much alteration.       

The rest of the paper is organised as follows. Section \ref{sec:problem} formally describes the problem and develops the optimization model. The formulations of the original Gini coefficient based on Lorenz curve, valid inequalities, and the cluster-based Lorenz curve are presented in Section \ref{subsec:equity}. Section~\ref{sec:results} presents the description of the Brazilian case study which is used to implement the model with real data and perform numerical analyses. The concluding section \ref{sec:conclusions} summarizes the paper's contributions and points out some remaining challenges and opportunities for future research.

\section{Problem set-up and mathematical formulation}
\label{sec:problem}

We propose a scenario-based two-stage stochastic programming model to optimize location-allocation decisions in disaster preparedness and response. Our problem concerns a given geographical area prone to natural hazards and is defined with respect to a number of settlements (e.g., a town or municipality, village or even neighbourhood), each of which we refer to as an disaster-prone area. In a disaster aftermath, these disaster-prone areas require relief aid from prepositioned stockpiles of supplies in response facilities (RFs) during the preparedness phase. For this purpose, humanitarian logisticians must decide, \textit{here-and-now}, on the sites to set-up RFs and the proper levels of critical supplies to maintain in stock at these RFs so as to effectively provide humanitarian assistance regardless of the potential disaster scenario. The \textit{wait-and-see}, or response, decisions reflect the assignment of victim needs, which we interchangeably term ``demands'', to the established RFs, given the occurrence of a disaster scenario. In a situation of scarce resources, wherein expenditures on preparedness and response activities are budgeted, the problem is to find an effective-equitable solution that not only only maximizes the total demand coverage, but also equitably apportions this coverage among disaster-prone areas.   

Mathematically, let $n \in N$ denote a potential location for an RF, $a \in A$ a disaster-prone area, $\ell \in L$ the size of an RF, $r\in R$ a type of relief aid, and $s\in S$ a disaster scenario. The parameters of the optimization model are as follows: An RF of size $\ell$ placed at location $n$ has storage capacity $\kappa^{\mbox{\scriptsize resp}}_{\ell n}$ (in volume) and an associated fixed cost $c^{\mbox{\scriptsize o}}_{\ell n}$. Each unit of relief aid $r$ requires $\rho_r$ units of storage space and incurs a unit prepositioning cost $c^{\mbox{\scriptsize p}}_{rn}$. The minimum prepositioning quantity of relief aid $r$ at RF location $n$ is given by $\theta^{\min}_{rn}$. There is a maximum quantity $\theta^{\max}_{r}$ of relief aid $r$ available for prepositioning. The unit cost of shipping relief aid supplies from RF location $n$ to disaster-prone area $a$ depends on the number of trips from $n$ to $a$ and is given by $c^{\mbox{\scriptsize d}}_{an}$. The vehicle that performs this activity has capacity $\kappa^{\mbox{\scriptsize v}}$. Parameters $\eta$ and $\eta'$ are the pre-disaster and post-disaster emergency relief funds available to carry out the logistics activities, respectively. The total victims' needs for relief aid $r$ at disaster-prone area $a$ in scenario $s\in S$ is given by $d_{ras}$, and the corresponding probability of occurrence of scenario $s$ is denoted by $\pi_s$, such that $\pi_s>0$ and $\sum_{s\in S}\pi_s=1$. We define the following first-stage decision variables: $Y_{\ell n}$ is a binary variable that indicates whether an RF of size $\ell$ is established at location $n$ $(Y_{\ell n}=1)$ or not $(Y_{\ell n}=0)$. $P_{rn}$ is the quantity of relief aid $r$ prepositioned at RF location $n$. The second-stage decision variable is $X_{rans}$, which represents the fraction of relief aid $r$ at disaster-prone area $a$ assigned to RF location $n$ in scenario $s$. We seek to generate decisions that maximize an {\it effectiveness $\times$ equity} objective, where our effectiveness function determines the extent to which the established RFs cover victim needs in a scenario $s$, and our equity function measures the extent to which the prepositioned stocks of relief aids are fairly allocated amongst disaster-prone areas in that scenario.  

The proposed optimization model is posed as follows:
\begin{alignat}{4}
\max\mbox{  } &  \Expec[Q({\bf Y},{\bf P},s)]\label{eq:priori}\\
\text{s.t. }\sum_{r\in R}\rho_r P_{rn} &\leq \sum_{\ell\in L}\kappa^{\mbox{\scriptsize resp}}_{\ell n} Y_{\ell n},\ \forall n\in N,\label{eq:locprep}\\
\sum_{n\in N } P_{rn} & \leq \theta^{\max}_{r},\ \forall r\in R,\label{eq:prep1}\\
\sum_{r \in R}P_{rn} & \geq \theta^{\min}_{n} \sum_{\ell\in L }Y_{\ell n},\ \forall r\in R,n\in N,\label{eq:prep2}\\
\sum_{\ell\in L } Y_{\ell n}&\leq 1,\ \forall n\in N, \label{eq:loc}\\
\sum_{r\in R,n\in N}c^{\mbox{\scriptsize p}}_{rn}P_{rn} + \sum_{\ell\in L,n\in N} c^{\mbox{\scriptsize o}}_{\ell n}Y_{\ell n} &\leq \eta,\label{eq:budget0}\\
Y_{\ell n} & \in \{0,1\},\ \forall \ell\in L,n\in N,\label{eq:domain1}\\
P_{rn} & \geq 0,\ \forall r\in R,n\in N,\label{eq:domain2}
\end{alignat}

\n where for every scenario $s\in S$,

\vspace{-0.5cm}

\begin{alignat}{4}
Q({\bf Y},{\bf P},s)&=\max\mbox{  } U_s(1-G_s)\label{eq:recourse}\\
\text{s.t. }\sum_{a\in A}X_{rans} d_{ras} &\leq P_{rn},\ \forall r\in R,n\in N,\label{eq:capacity1}\\
\sum_{n\in N}X_{rans} &\leq 1,\ \forall r\in R,a\in A,\label{eq:capacity2}\\
\sum_{r \in R,a\in A,n\in N}c^{\mbox{\scriptsize d}}_{an} \frac{\rho_r}{\kappa^{\mbox{\scriptsize v}}} d_{ras}X_{rans} &\leq \eta',\label{eq:budget1}\\
X_{rans}\ &  \geq 0,\ \forall r\in R,a\in A,n\in N\label{eq:domain4}.
\end{alignat}

The objective function~(\ref{eq:priori}) maximizes the expectation of the recourse function (there is no first-stage cost). Considering we have a discrete distribution of scenarios with finite support, $\Expec[Q({\bf Y},{\bf P},s)]=\sum_{s\in S}\pi_sQ({\bf Y},{\bf P},s)$. Constraint~(\ref{eq:locprep}) ensures that relief aid can be prepositioned at location $n$ only if an RF is established at that location. Constraint~(\ref{eq:prep1}) establishes a maximum quantity of relief aid $r$ to be prepositioned across all RFs. This quantity is generally defined a priori via agreements between private suppliers and public bodies or non-governmental organisations in charge of disaster relief operations. Constraint~(\ref{eq:prep2}) ensures that if prepositioning of relief aid $r$ takes place at RF location $n$, there is a minimum-quantity requirement. %This avoids sparsely stocked RFs. 
Constraint~(\ref{eq:loc}) ensures that each location $n$ can only admit an RF of a single size. Constraint~(\ref{eq:budget0}) defines the pre-disaster financial budget for RF setup and prepositioning. Constraints (\ref{eq:domain1}) and (\ref{eq:domain2}) specify the domains of the first-stage decision variables. The recourse function $ Q(\cdot)$ given in (\ref{eq:recourse}) trades-off effectiveness, via a demand coverage measure $U_s$, and equity, via the Gini coefficient $G_s$ (to be precise, the Gini coefficient is a measure of {\it inequity} and $(1-G_s)$ is a measure of equity). The exact forms of both metrics are discussed in the next section. Constraint (\ref{eq:capacity1}) ensures that the quantity of relief aid $r$ supplied from an RF at location $n$ does not exceed the prepositioned stock level at that RF. Constraint (\ref{eq:capacity2}) guarantees that the total quantity of relief aid $r$ supplied to disaster-prone area $a$ does not exceed victims' needs. Constraint (\ref{eq:budget1}) imposes a financial budget on post-disaster activities. Finally, constraint (\ref{eq:domain4}) specifies the domain of every second-stage decision variable.

\section{Equity modelling}
\label{subsec:equity}

\begin{figure}[H]
\centering
\includegraphics[scale=0.25]{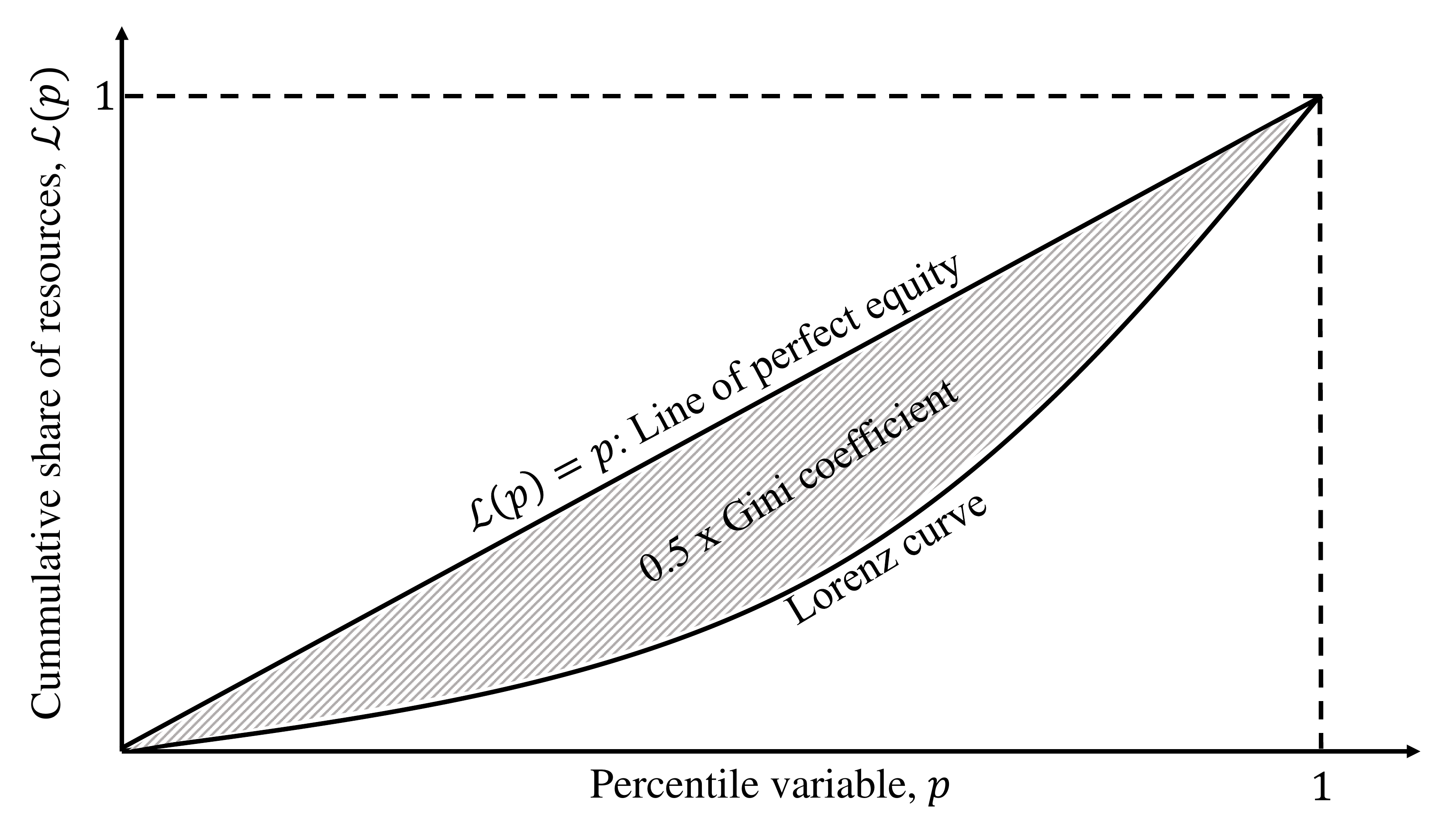}
\caption{Lorenz curve and the Gini coefficient.}
  \label{Lorenz}
\end{figure}

The basis for defining the Gini coefficient is the Lorenz curve, as shown in Figure \ref{Lorenz}. For consistency, we will use the definitions employed in \cite{farris2010gini}.  Suppose that a quantity of resources is distributed among groups in a population. If we sort these groups in increasing order of the shares of resources they possess, then $p$, the percentile variable, represents the poorest $100p\%$ of the population and $\mathcal{L}(p)$ is the cumulative share of resources they possess. If $\mathcal{L}(p)=p$ for all $p$, the distribution of resources is perfectly equitable. In reality, perfect equity rarely exists and the distribution of resources rather follows a convex, monotonically increasing, and scale-invariant curve that lies below the line of perfect equity, called the Lorenz curve. We note that while we assume convexity for the purposes of this work, strictly speaking, the Lorenz curve is not always convex, with nonconvexity occurring in restricted cases such as when resources within the population follow a very U-shaped Beta distribution \citep{prendergast2016}. The area between the line of perfect equity and the Lorenz curve is a measure of the inequity that exists in the population and the Gini coefficient is computed as twice this area. The reason for the factor-of-2 multiplication is to scale the area in such a way that the Gini coefficient is equal to 1 when one group in the population possesses all the resources, which represents the case of highest inequity (the area between the Lorenz curve and the line of perfect equity is then simply equal to the area of the triangle under the line of perfect equity, which is equal to 1/2). The mathematical definition of the Gini coefficient $G$ is 

\begin{alignat*}{4}
G \coloneqq 2\int_{0}^{1} [p- \mathcal{L}(p)]\,dp.
\end{alignat*}

\vspace{0.5cm}

In general, the integral cannot be directly evaluated because percentiles are only available for a {\it discrete} number of groups in the population, which means that the Lorenz curve is not fully defined over the range $[0,1]$. In order to fill the rest of the curve, linear interpolation is commonly employed, which is the simplest way to ensure continuity and integrability. 

The objective function of our model, $U_s(1-G_s)$, contains an effectiveness measure $U_s$ and a Gini coefficient $G_s$, and follows the rationale found in \cite{eisenhandler2019}. Unlike the Gini coefficient described in the previous paragraph, where the Lorenz curve is constructed from a priori available data, $G_s$ in our objective function is computed from optimal {\it decisions} that are not known a priori. The challenge we address in this paper is to translate the underlying steps of Gini coefficient computation, i.e. the construction of the percentile variable and its associated linearly interpolated Lorenz curve and the evaluation of the area between the Lorenz curve and the line of perfect equity, to a decision-driven setting in such a way that our optimization model remains solvable. Before addressing this challenge, we need to establish the metric which we aim to `equitabilize', starting from our effectiveness measure, which we define as  

\begin{alignat*}{4}
U_s\coloneqq\sum_{r \in R,a\in A,n\in N}u_{ras}X_{rans}, \ \forall s\in S,
\end{alignat*}

\vspace{0.5cm}

\noindent in which $u_{ras}=d_{ras}/\sum_{r' \in R,a' \in A}d_{r'a's},\ \forall r\in R, a \in A, s \in S$. As such, $u_{ras}$ represents the fraction of the total demand that is for relief aid $r$ in disaster-prone area $a$ in scenario $s$, and $U_s$ is the fraction of the total demand that is served/covered in scenario $s$. The quantity $\sum_{r \in R,n\in N}u_{ras}X_{rans}$ is thus the fraction of the total demand that is covered in disaster-prone area $a$ in scenario $s$, which we term the demand coverage at disaster-prone area $a$. Our aim is to provide equitable demand coverage among disaster-prone areas, without neglecting the effectiveness measure, which is achieved with the maximization of $U_s(1-G_s)$. Indeed, this objective function avoids, among others, the case with zero demand coverage in all areas, which despite being perfectly equitable, is completely ineffective. 

To generate our Lorenz curve for every scenario $s \in S$ in a decision-driven setting, we first start by deriving percentile values with decision variables. If disaster-prone area $a$ has the $j^{th}$ ranked demand coverage among the disaster-prone areas in set $A_s$, where $A_s=\{a\in A: \sum_{r\in R}d_{ras} > 0\}$, its percentile value is $p_{as}=j/|A_s|$, where $|A_s|$ is the cardinality of set $A_s$. Letting $O_{ajs}$ be a binary decision variable that takes a value of $1$ if disaster-prone area $a$ has the $j^{th}$ ranked demand coverage and a value of $0$ otherwise, we can rewrite $p_{as}$ as

\begin{alignat*}{4}
&p_{as} = \dfrac{\sum_{j\in[|A_s|]} jO_{a js}}{|A_s|} ,\ \forall a \in A_s,\\
\text{s.t. }&\sum_{j\in[|A_s|]} O_{ajs} = 1,\ \forall a \in A_s,\label{P1}\\
&\sum_{a\in A_s} O_{ajs} = 1,\ \forall j\in[|A_s|], \\
&\sum_{r\in R,n\in N} u_{ras} X_{rans}O_{ajs} \le \sum_{r\in R,n\in N} u_{ras} X_{rans}O_{aj+1s} ,\ \forall a \in A_s, j\in[|A_s|-1], \\
&O_{ajs} \in\{0,1\},\ \forall a \in A_s, j\in [|A_s|], 
\end{alignat*}

\vspace{0.5cm}

\noindent where notation $[n]$ represents the set of running indices from 1 to $n$. This means that the formula $(\sum_{j\in[|A_s|]} jO_{a js})/|A_s|$ computes the percentile value of disaster-prone area $a$ if the binary variable is constrained such that \textit{i}) every disaster-prone area with positive total demand is assigned to a unique rank, \textit{ii}) every rank is assigned to a unique disaster-prone area with positive total demand, and \textit{iii}) the $j^{th}$ ranked demand coverage does not exceed the $(j+1)^{th}$ ranked demand coverage. Constraint $\sum_{r\in R,n\in N} u_{ras} X_{rans}O_{ajs} \le \sum_{r\in R,n\in N} u_{ras} X_{rans}O_{aj+1s}$, $\forall a \in A_s, j\in[|A_s|-1]$, is nonlinear since it involves the multiplication of two decision variables. An exact linearization is possible by defining a new decision variable $Z_{js}$ that represents the demand coverage of the $j^{th}$ ranked area and replacing the nonlinear constraints with the list of constraints 

\begin{alignat*}{4}
&Z_{js} \le \sum_{r\in R,n\in N}u_{ras}X_{rans} + 1 - O_{ajs},\ \forall a\in A_s,j\in [|A_s|], \\
&Z_{js} \ge \sum_{r\in R,n\in N}u_{ras}X_{rans} - 1 + O_{ajs},\ \forall a\in A_s,j\in [|A_s|], \\
&Z_{js} \le Z_{j+1 s},\ \forall a\in A_s,j\in [|A_s|-1], \\
&Z_{js} \ge 0,\ \forall j\in [|A_s|].
\end{alignat*}

\vspace{0.5cm}

We can see from these constraints that if disaster-prone area $a$ has the $j^{th}$ ranked demand coverage, $O_{ajs}=1$, which means that $Z_{js}\le \sum_{r\in R,n\in N}u_{ras}X_{rans}$ and $Z_{js}\ge \sum_{r\in R,n\in N}u_{ras}X_{rans}$, which is satisfied if and only if $Z_{js}= \sum_{r\in R,n\in N}u_{ras}X_{rans}$. On the other hand, given that $0\le Z_{js}\le 1$ by definition, if $O_{ajs}=0$, the constraints do not restrict $Z_{js}$ any way, since $\sum_{r\in R,n\in N}u_{ras}X_{rans} - 1 \le 0$  and $\sum_{r\in R,n\in N}u_{ras}X_{rans} + 1 \ge 1$.

\begin{figure}[H]
\centering
\includegraphics[scale=0.25]{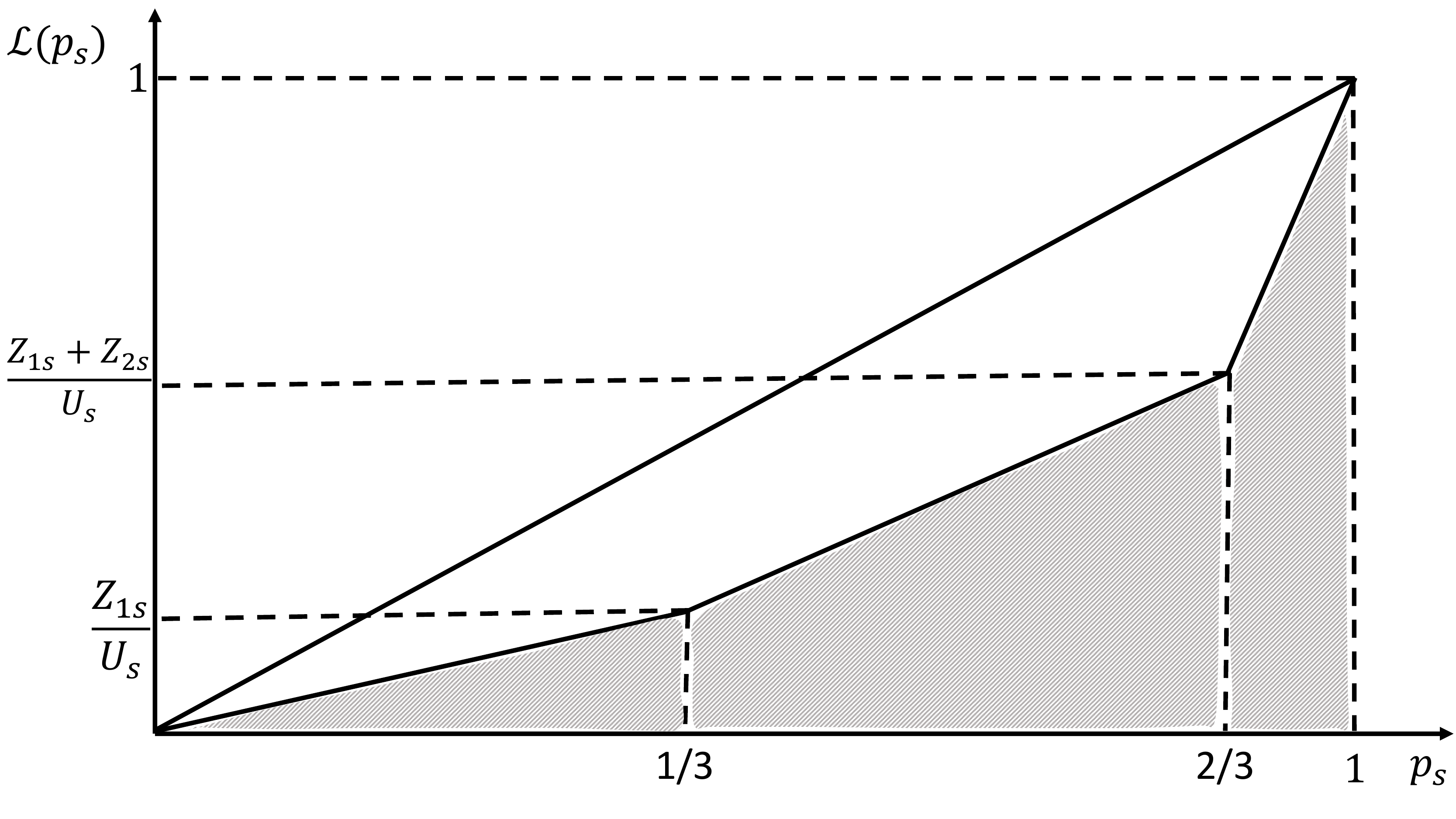}
\caption{Decision-driven linearly interpolated Lorenz curve and its shaded adjacent trapezoids.}
  \label{lin}
\end{figure}

With the above linearized representation of percentile values, the cumulative share of demand coverages up to the $j^{th}$ rank is $(\sum_{j'\in [j]}Z_{j's})/U_s$, and the resulting decision-driven linearly interpolated Lorenz curve form adjacent trapezoids with the percentile value axis, as illustrated in Figure \ref{lin} for a three-disaster-prone-area setting. The area between the line of perfect equity and the Lorenz curve (which is half the Gini coefficient) is then the difference between the area of the triangle under the line of perfect equity $(=1/2)$ and the total area formed by the trapezoids.

The area under the Lorenz curve between the $(j-1)^{th}$ and the $j^{th}$ percentile value under linear interpolation, for $j\in [|A_s|]$, where $\mathcal{L}(0)=0$, is the area of a trapezoid, which is

\begin{alignat*}{4}
&\dfrac{Z_{1s}}{2 |A_s| U_s}  \text{, if }  j=1\\
&\dfrac{\sum_{j^{'}\in[j-1]}  Z_{j's} + \sum_{j'\in [j]}  Z_{j's}}{2 |A_s| U_s}  ,\ \forall j \in [|A_s|]\setminus \{1\}.
\end{alignat*}

\vspace{0.5cm}

The Gini coefficient is therefore

\begin{alignat*}{4}
G_s=1 - \dfrac{1}{|A_s| U_s}\big[Z_{1s}+\sum_{j\in [|A_s|]\setminus \{1\}}(\sum_{j'\in[j-1]}  Z_{j's} + \sum_{j'\in [j]}  Z_{j's})\big], 
\end{alignat*}

\vspace{0.5cm}

\noindent which implies that 

\begin{alignat*}{4}
& U_{s} (1-G_s)&&=\dfrac{1}{|A_s|}\big[Z_{1s}+\sum_{j\in [|A_s|]\setminus \{1\}}(\sum_{j'\in[j-1]}  Z_{j's} + \sum_{j'\in [j]}  Z_{j's})\big]\\
& && =\sum_{j \in [|A_s|]}\dfrac{1}{|A_s|}(2|A_s| +1- 2j) Z_{js}.
\end{alignat*}

\vspace{0.5cm}

Finally, we formulate the deterministic equivalent of our original two-stage stochastic programming model as 

\begin{alignat}{4}
\max\mbox{  }   \sum_{j \in [|A_s|],s\in S}&\dfrac{1}{|A_s|}(2|A_s| +1- 2j) \pi_s Z_{js}\\
\text{s.t. } \eqref{eq:locprep} &- \eqref{eq:domain2},\nonumber\\ 
\eqref{eq:capacity1}& - \eqref{eq:domain4},\ \forall s\in S,\nonumber\\
\sum_{j\in[|A_s|]} O_{ajs} &= 1,\ \forall a \in A_s,s\in S,\\
\sum_{a\in A_s} O_{ajs} &= 1,\ \forall j\in[|A_s|],s\in S, \\
Z_{js} &\le \sum_{r\in R,n\in N}u_{ras}X_{rans} + 1 - O_{ajs},\ \forall a\in A_s,j\in [|A_s|],s\in S, \\
Z_{js} &\ge \sum_{r\in R,n\in N}u_{ras}X_{rans} - 1 + O_{ajs},\ \forall a\in A_s,j\in [|A_s|],s\in S, \\
Z_{js} &\le Z_{j+1 s},\ \forall j\in [|A_s|-1],s\in S, \\
O_{ajs} &\in\{0,1\},\ \forall a \in A_s, j\in [|A_s|],s\in S,\\
Z_{js} &\ge 0,\ \forall j\in [|A_s|],s\in S.
\end{alignat}

\subsection{Valid inequality based on an upper-bounding Lorenz curve}
The Gini coefficient calculation relies on ranking demand coverages of disaster-prone areas. Because no lower bounds are specified on allocation variables in our model, there are no clear dominance relationships among the rankings. This can lead to loose linear relaxations and severely suboptimal root nodes in solution algorithms, which generally lead to longer solution times and lower numerical efficiency. Here, we develop a valid inequality, through an upper-bounding Lorenz curve, to eliminate some dominated area rankings and thus tighten the linear relaxation bound. 

Let $H^L_s \ge 0$ be a decision variable whose maximum value is the $1^{st}$ ranked (lowest) demand coverage and $H^{U}_s \ge 0$ be a decision variable whose minimum value is the $|A_s|^{th}$ ranked (highest) demand coverage. By definition,

\begin{alignat*}{4}
&H^{L}_s &&\le \sum_{r \in R,n\in N} u_{ras}X_{rans},\ \forall a \in A_s,\  s\in S,\\
&H^{U}_s &&\ge \sum_{r \in R,n\in N} u_{ras}X_{rans},\ \forall a \in A_s,\ s\in S.
\end{alignat*}  

\vspace{0.5cm}

\begin{figure}[H]
\centering
\includegraphics[scale=0.25]{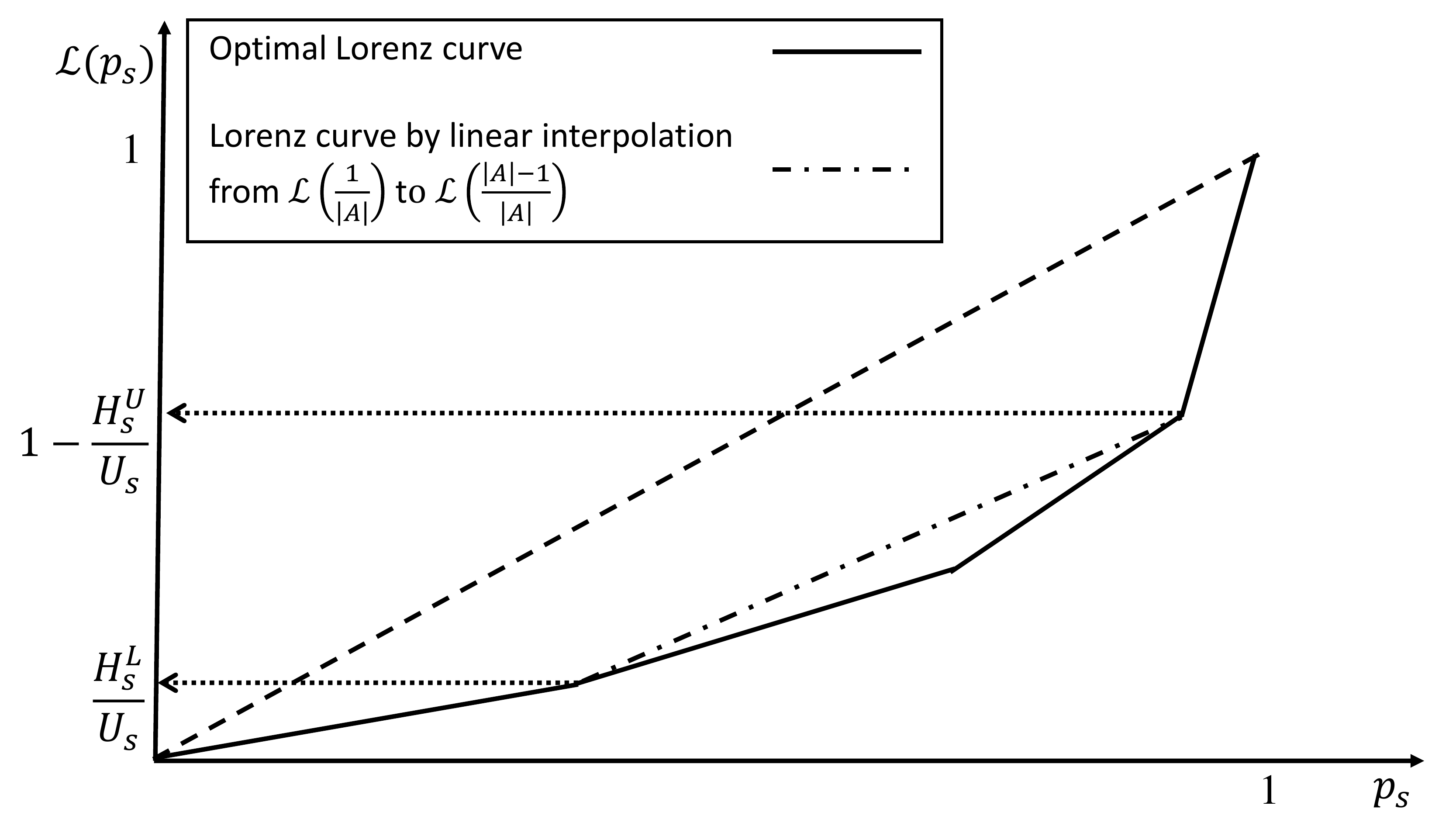}
\caption{Illustration of upper-bounding Lorenz curve.}
  \label{Gini_Bound1}
\end{figure}

Suppose that our deterministic equivalent with the above two constraints yields an optimal disaster-prone area ranking, and therefore an optimal Lorenz curve, such that $\mathcal{L}(\frac{1}{|A_s|})=\frac{H^{L}_s}{U_s}$ and $\mathcal{L}(\frac{|A_s|-1}{|A_s|})=1-\frac{H^{U}_s}{U_s}$. It is clear that by linear interpolation from $\mathcal{L}(\frac{1}{|A_s|})$ to $\mathcal{L}(\frac{|A_s|-1}{|A_s|})$, we obtain an upper-bounding Lorenz curve to the optimal one, as shown in Figure \ref{Gini_Bound1}.

The sum of the areas of the 3 trapezoids underneath the upper-bounding Lorenz curve is 

\begin{alignat*}{4}
&\dfrac{1}{2U_s|A_s|}H^L_s + \dfrac{|A_s|-2}{2U_s|A_s|}(H^L_s + U_s - H^U_s) + \dfrac{1}{2U_s|A_s|}(U_s-H^U_s+U_s),
\end{alignat*}

\vspace{0.5cm}

\noindent which gives 

\begin{alignat*}{4}
&U_{s} (1-G_s) &&=\dfrac{1}{|A_s|}\big(H^L_s + (|A_s|-2)(H^L_s + U_s - H^U_s) + 2U_s - H^U_s \big)\\
& &&= \dfrac{1}{|A_s|}\big(|A_s|U_s+(|A_s|-1)(H^L_s- H^U_s)  \big),
\end{alignat*}

\vspace{0.5cm}

\noindent from which, because the upper-bounding Lorenz curve has smaller Gini coefficient for a fixed effectiveness measure, we devise the valid inequality

\vspace{0.5cm}

\begin{alignat}{4}
&\sum_{j\in [|A_s|]}(2|A_s| +1- 2j) Z_{js} &&\le |A_s|U_s + (|A_s|-1)(H^L_s-H^U_s). \label{valI}
\end{alignat}

\subsection{Cluster-based Lorenz curve}
One of the main difficulties in building the decision-driven Lorenz curve is that it may require an exploration of $|A_s|!$ ranks of disaster-prone areas, with this process being further compounded by the ranks being decision-dependent. In this section, we propose to use a priori $k_s$-means clustering and construct our decision-driven Lorenz curve based on clusters, rather than individual disaster-prone areas. The main underlying idea here is that disaster-prone areas with high degrees of similarity are likely to be close in percentile to each other, and thus the resulting cluster-based Lorenz curve will likely be an upper-bounding Lorenz curve of the optimal one. If the cluster-based Lorenz curve is close enough to the optimal one, we expect the optimal decisions from the two models (with optimal Lorenz curve and with cluster-based Lorenz curve)  to not differ significantly. In addition, if $A_s=A$ (all disaster-prone areas have positive total demand) and $k_s$-means clustering is performed such that $k_s=|A|$, $\forall s \in S$, every resulting cluster will contain a single and unique disaster-prone area, which means that the cluster-based Lorenz curve will be exactly the same as the optimal one. We note that the main purpose here is not to offer a theoretical treatise of this idea, but rather to use it to make the two-stage stochastic programming model with Lorenz-curve-based Gini more translatable to practice. Moreover, our empirical tests conducted in the following section report promising results, both in terms of solution times and equity, compared to the benchmark Gini mean difference approach that is traditionally used in the literature to circumvent tractability issues. Our cluster-based approach is data-driven and improves numerical efficiency without sacrificing the fundamental principles used in constructing the Lorenz curve, whereas the Gini mean difference approach entirely avoids constructing the Lorenz curve in order to improve solvability.  

Let $B_s = \{B_{1s},B_{2s},\dots,B_{k_ss}\}$, be an ordered set of clusters, where the clusters are obtained from the $k_s$-means clustering approach, and $B_{ws} \subseteq A_s$, $|B_{ws}|\ge  1$,  $\forall w \in [k_s]$, $\bigcup_{w \in [k_s]}B_{ws} = A_s$ and $B_{ws}\bigcap B_{w's} = \{ \}$, $\forall w,w' \in [k_s], w\neq w'$. Our deterministic equivalent under cluster-based Lorenz curve is derived in the same way as our original deterministic equivalent, except that the Lorenz curve is constructed on clusters of disaster-prone areas, rather than individual ones. For brevity, we will not repeat the derivation here. The deterministic equivalent, hereafter called Gini with Clusters (GiniC), is then cast as 

\begin{alignat}{4}
\max\mbox{  }  \sum_{j \in [k_s],s\in S}&\dfrac{1}{k_s}(2k_s +1- 2j) \pi_s Z_{js}\\
\text{s.t. } \eqref{eq:locprep} &- \eqref{eq:domain2}, \nonumber\\
\eqref{eq:capacity1}& - \eqref{eq:domain4},\ \forall s\in S,\nonumber\\
\sum_{j\in[k_s]} O_{wjs} &= 1,\ \forall w \in [k_s],s\in S,\\
\sum_{w\in [k_s]} O_{wjs} &= 1,\ \forall j\in[k_s],s\in S, \\
Z_{js} &\le \sum_{r\in R,n\in N,a\in B_{ws}}u_{ras}X_{rans} + 1 - O_{wjs},\ \forall w\in [k_s],j\in [k_s],s\in S, \\
Z_{js} &\ge \sum_{r\in R,n\in N,a\in B_{ws}}u_{ras}X_{rans} - 1 + O_{wjs},\ \forall w\in [k_s],j\in [k_s],s\in S, \\
Z_{js} &\le Z_{j+1 s},\ \forall j\in [k_s-1],s\in S, \\
O_{bjs} &\in\{0,1\},\ \forall b \in [k_s], j\in [k_s],s\in S,\\
Z_{js} &\ge 0,\ \forall j\in [k_s],s\in S,
\end{alignat}

\vspace{.5cm}

\noindent where $O_{wjs}=1$ if the $w^{th}$ cluster has the $j^{th}$ ranked demand coverage and $O_{wjs}=0$ otherwise, and the demand coverage for the $w^{th}$ cluster is calculated as $\sum_{r\in R,n\in N,a\in B_{ws}}u_{ras}X_{rans}$. Since ranking is cluster-based, this model explores $k_s!$ ranks instead of $|A_s|!$ ranks and when $k_s < |A_s|$, for some $s \in S$, the number of ranks decreases fast. Furthermore, since $k_s$ is user-specified, the decision maker can control the trade-off between the numerical efficiency of the resulting model and the accuracy with which the optimal Lorenz curve is approximated with the cluster-based Lorenz curve.

\section{Results}
\label{sec:results}

The proposed models were coded in GAMS 25.1.1 software and solved with CPLEX 12.8 (default settings) on a computer with 16GB RAM, Intel Core i7 and Windows 7 operating system. The stopping criteria are either elapsed times exceeding 3600 seconds or optimality gaps relative to the best lower bound smaller than 0.001\%. We provide results for an application in Brazil and show the benefits of our approach over existing methods.   

\subsection{Application in Brazil}
\label{sec:application}

We are particularly motivated by the challenging problem of allocating scarce disaster relief aid to vulnerable areas in \emph{Serrana} region of Rio de Janeiro state in Brazil. We choose the Serrana Region for our case-study for several reasons: it is a mountainous area prone to natural hazards due its geomorphology and climate; environmental degradation combined with unplanned land use and occupation caused by socioeconomic problems make its residents vulnerable to natural hazards; torrential rains followed by floods and landslides/mudslides are frequently recurring events that strike this area and affect thousands of people every year \citep{Moreno2016,AlemClarkMoreno2016}; the \emph{Serrana} region is the Brazilian area with the highest number of fatal victims caused by natural hazards \citep{Coppetec2014,s2id2020}; the so-called Megadisaster of the \emph{Serrana} Region in January 2011 caused hundreds of fatal victims and missing people, as well as thousands of homeless and displaced people and is considered the worst disaster ever recorded in Brazil in terms of fatal victims \citep{CameraRJ}; it is worth mentioning that this Megadisaster is also among the ten worst landslides worldwide since 1900.

Our case-study data consists of 13 disaster-prone areas (municipalities): Areal (\textsc{are}), Bom Jardim (\textsc{bjd}), Cordeiro (\textsc{cor}), Macuco (\textsc{mac}), Nova Friburgo (\textsc{nfb}), Petrópolis (\textsc{pet}), Santa Maria Madalena (\textsc{smm}), São José do Vale do Rio Preto (\textsc{srp}), São Sebastião do Alto (\textsc{ssa}), Sapucaia (\textsc{sap}), Sumidouro (\textsc{sum}), Teresópolis (\textsc{ter}), and Três Rios (\textsc{trr}). We use the historical data from the period 2000--2018 to estimate victim needs, building a total of 18 equiprobable scenarios, as shown in Table~\ref{tab:sto_victims} of \ref{ap:input_data}. This information was obtained from the Integrated Disaster Information System \citep{s2id2020}, a platform that provides the records of the National Secretariat for Civil Protection and Defense (SEDEC) under the Brazilian Ministry of Regional Development.

We consider that there is a basic requirement for six types of relief aid: water, food, hygiene kits, cleaning kits, dormitory kits, and mattress. Each aid has a \textit{coverage} in terms of number of people, e.g., a kit of food covers the nutritional requirements of a four-person family. The demand $d_{ras}$, which represents the victim needs in terms of the quantity of relief aid $r$ for each disaster-prone area $a$ in scenario $s$, is evaluated as $
d_{ras} =  \left\lceil \frac{\mbox{length}_r}{\mbox{coverage}_r}\times \mbox{victims}'_{as}\right\rceil
$, where `$\mbox{length}_r$' is the number of days in which victims need to be supplied with relief aid $r$; $\mbox{victims}'_{as}$ is the number of homeless and displaced people in municipality $a$ in scenario $s$ shown in Table~\ref{tab:sto_victims}; $\mbox{coverage}_r$ is the number of people covered by one unit of relief aid $r$. We use the information about relief aid items referred to in the public announcement for procurement of humanitarian supplies by the Brazilian government \citep{ata2017} and used in recent studies in the literature \citep{Alem2020}. All the characteristics of the relief aid are summarised in Table~\ref{sec:aid} of Appendix A.

We consider four possible sizes for the response facilities in the different locations: small, medium, large and very large, whose capacities $\kappa^{\mbox{\scriptsize resp}}_{\ell n}$ and associated fixed cost $c^{\mbox{\scriptsize o}}_{\ell n}$ are shown in Table~\ref{tab:facilities} in \ref{ap:input_data}. The cost $c^{\mbox{\scriptsize o}}_{\ell n}$ for establishing a response facility of size $\ell$ was assumed proportional to its construction cost. The minimum prepositioning quantity was set to 1.

Shipping costs ($c^{\mbox{\scriptsize d}}_{an}$) were evaluated based on the assumption that relief aids are transported by medium-sized trucks, each with a capacity of $12m^3$. These vehicles mainly use  diesel as fuel, the cost of which is BRL 3.59 per litre \citep{anp2018}. We also assume that the trucks have an average consumption rate of 2.5km/litre. Finally, the unit shipping cost is calculated as $c^{\mbox{\scriptsize d}}_{an}=\frac{\mbox{diesel cost}}{\mbox{consumption}} \times \mbox{dist}_{an}$, where $\mbox{dist}_{an}$ is the distance (in km) between two nodes $a-n$ obtained using the \textit{Openrouteservice} (\texttt{https://openrouteservice.org/}). This service is useful to compute many-to-many distances and is based on data from the \textit{OpenStreetMap}, an open initiative to create and provide free geographic data. In the absence of the exact address of a response facility, the distances were approximated by using the centroid of each area. Table \ref{table:distance} in \ref{ap:input_data} shows all the pairwise distances.

We consider financial budgets of BRL 26,206,190 and BRL 23,415 for first- and second-stage decisions, respectively. These were calculated by solving a minimum-cost model such that approximately 30\% of victims' needs are satisfied, which is aligned with situations with very scarce resources.

For the GiniC model (model with cluster-based Lorenz curve), $k_s$-means clustering is performed for each scenario on a dataset comprising $\sum_{r\in R}u_{ras}$, $\forall a \in A$. The appropriate value of $k_s$, $\forall s \in S$, is established via identifying the elbow in the clustering scree plot. This gives the result $(k_1,k_2,\dots,k_{18})=(1,1,1,2,2,3,1,3,3,3,3,3,3,2,1,2,2,3).$ 

\subsection{Benchmark methods}
\label{sec:comparison}

We compare our models to two benchmark approaches: (i) a stochastic programming model without equity (hereafter simply called Stochastic Problem or SP), whose goal is solely the maximization of the effectiveness; and (ii) a stochastic programming model where the equity measure is the popular Gini Mean Difference formulation (hereafter called GMD). While several variants of the Gini mean difference formulation exist. e.g., \citet{mandell1991,eisenhandler2019,Mostajabdaveh2019}, we will use the most popular of them, which is developed by \cite{mandell1991}. The formulations of these benchmark models are provided below.

\begin{itemize}

\item SP is formulated as
\begin{align*}
\max_{} &\sum_{s\in S}\pi_sU_s\\ \text{s.t. } \eqref{eq:locprep} &- \eqref{eq:domain2}, \nonumber\\
\eqref{eq:capacity1}& - \eqref{eq:domain4},\ \forall s\in S.\nonumber
\end{align*}

\item GMD is cast as 

\begin{alignat*}{4}
\max  &\mbox{ } \sum_{s\in S}\pi_s\big(\sum_{r\in R,a \in A,n\in N}  u_{ras}X_{rans} - \sum_{a \in A,a' \in \{A:a'>a\}} t_{aa's} \big)\\
\text{s.t. } \eqref{eq:locprep} &- \eqref{eq:domain2}, \nonumber\\
\eqref{eq:capacity1}& - \eqref{eq:domain4},\ \forall s\in S,\nonumber\\
t_{aa's} &\ge\rho_{as}\sum_{r\in R,n\in N}  u_{ra's}X_{ra'ns} - \rho_{a's}\sum_{r\in R,n\in N} u_{ras} X_{rans},\,\forall a,a'\in A, s \in S,\\
t_{aa's} &\ge  \rho_{a's}\sum_{r\in R,n\in N} u_{ras} X_{rans} -  \rho_{as}\sum_{r\in R,n\in N} u_{ra's}X_{ra'ns}, \,\forall a,a'\in A, s \in S,\\
\end{alignat*}
where $\rho_{as}=\big(\sum_{r\in R}u_{ras}\big)/\big(\sum_{'r\in R,a' \in A}u_{r'a's}\big)$ for disaster-prone area $a$ in scenario $s$ is what is called the ``proportion of equity units available''. The \citet{mandell1991} formulation of Gini for our case would be
$$G_s=\dfrac{1}{U_s}\sum_{a \in A,a' \in \{A:a'>a\}}|\rho_{a's}\sum_{r\in R,n\in N}  u_{ras}X_{rans}-\rho_{a's}\sum_{r\in R,n\in N}  u_{ra's}X_{ra'ns}|,$$ 

which in linearized form, yields Model GMD. One can clearly see the better tractability of the GMD model, which does not contain additional binary variables and requires no ranking procedure, compared to the model with the original Gini coefficient based on the Lorenz curve.

\end{itemize}

\subsection{Analysis of the solutions}
\label{sec:solution}

\begin{figure}[H]
\centering
\includegraphics[scale=0.30]{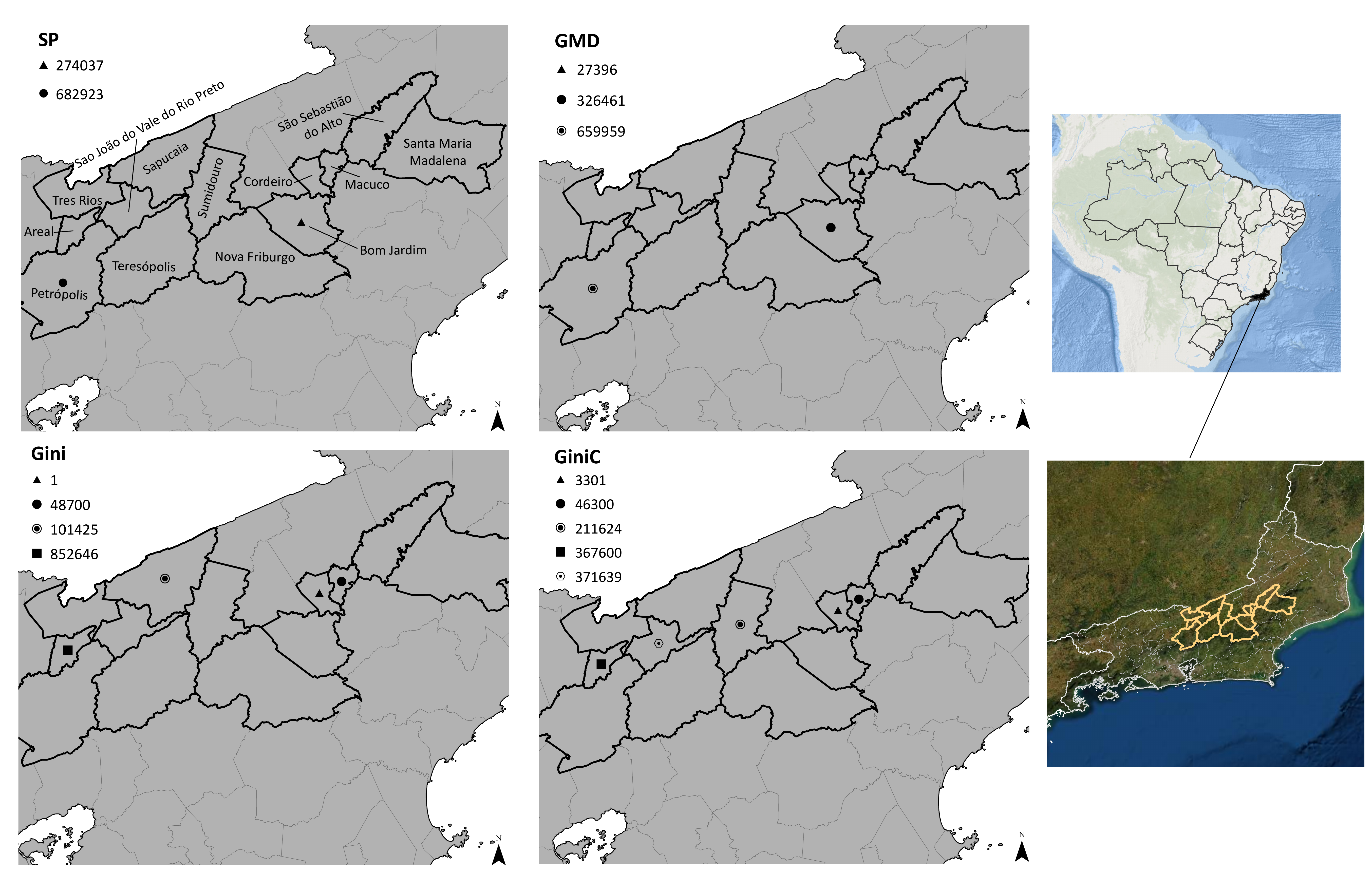}
\vspace{-0.5cm}
\caption{Map showing the Serrana region of Rio de Janeiro state. The optimal first-stage solution (location of RFs and prepositioning of relief aid) given by each approach is also indicated.}
  \label{fig:map}
\end{figure}

 Figure \ref{fig:map} shows the optimal setup and aggregate prepositioning decisions for all the proposed approaches. To simplify discussions, from now on, we will refer to our model with original Gini coefficient formulation via the Lorenz curve as ``Gini''. Overall, the quantity of prepositioned relief aid is quite similar among the four models ($\approx$ 1M units of items). However, the location of the response facilities and stockpiles varies greatly. Whereas SP establishes only two RFs (Bom Jardim and Petropolis), GiniC locates five RFs (Areal, Cordeiro, Macuco, Sao Jose do Vale do Rio Preto, and Sumidouro), which suggests that the decentralization of RFs help to achieve better equity levels. 
 
Table \ref{tab:stats} compares the four models in terms of coverage. We see that Gini and GiniC perform much better compared with SP and GMD approaches. In particular, GiniC manages to improve the worst-case coverage by almost four times when compared to GMD. Even though GMD provides an average  coverage slightly superior than SP ($\approx$ 2\%), the perfect coverage performance of the latter is 22\% better. The poor performance of GMD in terms of perfect coverage is notably evident for mattress (9.46\% versus 29.7\% given by the best approach). It is also worth noting that coverage levels given by Gini and GiniC are more equitable among relief aid items than those produced by either SP or GMD. The coefficient of variation (CoV) confirms this result (19.9\% given by GiniC versus 76.7\% given by GMD). 

\begin{table}[H]
\footnotesize
  \centering
    \caption{Coverage values and statistics for the proposed approaches.}
    \begin{tabular}{lrrrrrrrr}\midrule
          & \multicolumn{4}{c}{Coverage$^a$} & \multicolumn{4}{c}{Perfect coverage$^b$} \\
    \midrule
          & \multicolumn{1}{c}{SP} & \multicolumn{1}{c}{GMD} & \multicolumn{1}{c}{Gini} & \multicolumn{1}{c}{GiniC} & \multicolumn{1}{c}{SP} & \multicolumn{1}{c}{GMD} & \multicolumn{1}{c}{Gini} & \multicolumn{1}{c}{GiniC} \\
    \midrule
    Food  & \textit{0.1827} & 0.2486 & 0.5189 & \textbf{0.5420} & \textit{0.1351} & 0.1757 & \textbf{0.4730} & \textbf{0.4730} \\
    Water & \textit{0.7873} & 0.8332 & 0.8602 & \textbf{0.8816} & 0.7297 & \textit{0.5405} & 0.7297 & \textbf{0.7973} \\
    Hygiene & \textit{0.2420} & 0.2463 & 0.5260 & \textbf{0.5504} & 0.2027 & \textit{0.1892} & \textbf{0.5000} & 0.4865 \\
    Cleaning  & 0.2813 & \textit{0.2373} & 0.5562 & \textbf{0.6109} & 0.2297 & \textit{0.1757} & 0.5270 & \textbf{0.5676} \\
    Mattress & 0.1935 & \textit{0.1368} & 0.3630 & \textbf{0.6109} & 0.1622 & \textit{0.0946} & \textbf{0.2973} & 0.2703 \\
    Dormitory & \textit{0.2404} & 0.2665 & 0.5633 & \textbf{0.5974} & 0.1757 & \textit{0.1622} & 0.5270 & \textbf{0.5541} \\
    \midrule
    Average & 0.3212 & 0.3281 & 0.5646 & 0.6322 & 0.2725 & 0.2230 & 0.5090 & 0.5248 \\
    St. dev. & 0.2312 & 0.2517 & 0.1623 & 0.1258 & 0.2264 & 0.1591 & 0.1383 & 0.1708 \\
    CoV$^c$ (\%)  & 71.96 & 76.71 & 28.74 & 19.90 & 83.06 & 71.38 & 27.16 & 32.55 \\
    Best-case   & 0.7873 & 0.8332 & 0.8602 & 0.8816 & 0.7297 & 0.5405 & 0.7297 & 0.7973 \\
    Worst-case   & 0.1827 & 0.1368 & 0.3630 & 0.5420 & 0.1351 & 0.0946 & 0.2973 & 0.2703 \\
    \bottomrule
    \end{tabular}%    	
    \begin{tablenotes}
		\item $^a$ `Coverage' represents the average of coverage across disaster-prone areas and scenarios.
		\item $^b$ `Perfect coverage' shows the percent of times in which the coverage is 100\%.
		\item $^c$ `Coefficient of Variation' is measured as the ratio (in \%) between standard deviation and average values.
		Note. \textbf{Best} results. \textit{Worst} results.
	\end{tablenotes}
  \label{tab:stats}%
\end{table}%

\subsection{Effectiveness versus equity analysis}
\label{sec:tradeoff}

The analysis of the trade-offs between effectiveness and equity of the proposed approaches is performed by means of the following procedure. First, we solve the problem (SP/GMD/Gini/GiniC) considering the original set of scenarios $s \in S$. Then, we randomly generate 100 realizations of victim needs following a uniform distribution in the interval $\displaystyle \big[ \min_{s \in S} \{d_{ras}\}, \max_{s \in S} \{d_{ras}\} \big]$, $\forall r \in R$ and $a \in A$. We fix the first-stage decision variables ({\bf Y},{\bf P}) obtained with the original scenarios, and evaluate the problem for each random realization $s'$. Finally, the Gini coefficient and the effectiveness are evaluated as follows: 
\begin{enumerate}
    \item Take the optimal solution $X^*_{rans'}$ and compute $X^{\text{rank}}_{as'} = \sum_{r \in R,n \in N} u_{ras'}X^*_{rans'}$, where $\displaystyle u_{ras'} = \frac{d_{ras'}}{\sum_{r \in R,a \in A}d_{ras'}}$.
    \item Rank $X^{\text{rank}}_{as'}$ over $a\in A_{s'}$ in ascending order and store the ranked sequence as $(Z^*_{1s'}, Z^*_{2s'}, \dots, Z^*_{|A_{s'}|s'})$, where $Z^*_{js'}$ is the $j^{th}$ ranked value. 
    \item Compute the optimal Gini coefficient as 
    \begin{align}
        G^*_{s'}=1 - \dfrac{1}{|A_{s'}|\sum_{a \in A_{s'}}X^{\text{rank}}_{as'}}\big[Z^*_{1s'}+\sum_{j\in [|A_{s'}|]\setminus \{1\}}(\sum_{j'\in[j-1]}  Z^*_{j's'} + \sum_{j'\in [j]}  Z^*_{j's'})\big]\label{eq:G}
    \end{align}
      and the effectiveness as 
      \begin{align}
      U^*_{s'}=\sum_{a \in A}X^{\text{rank}}_{as'}\label{eq:U}.
        \end{align}
\end{enumerate}

\begin{figure}[H]
\centering
\includegraphics[scale=0.45]{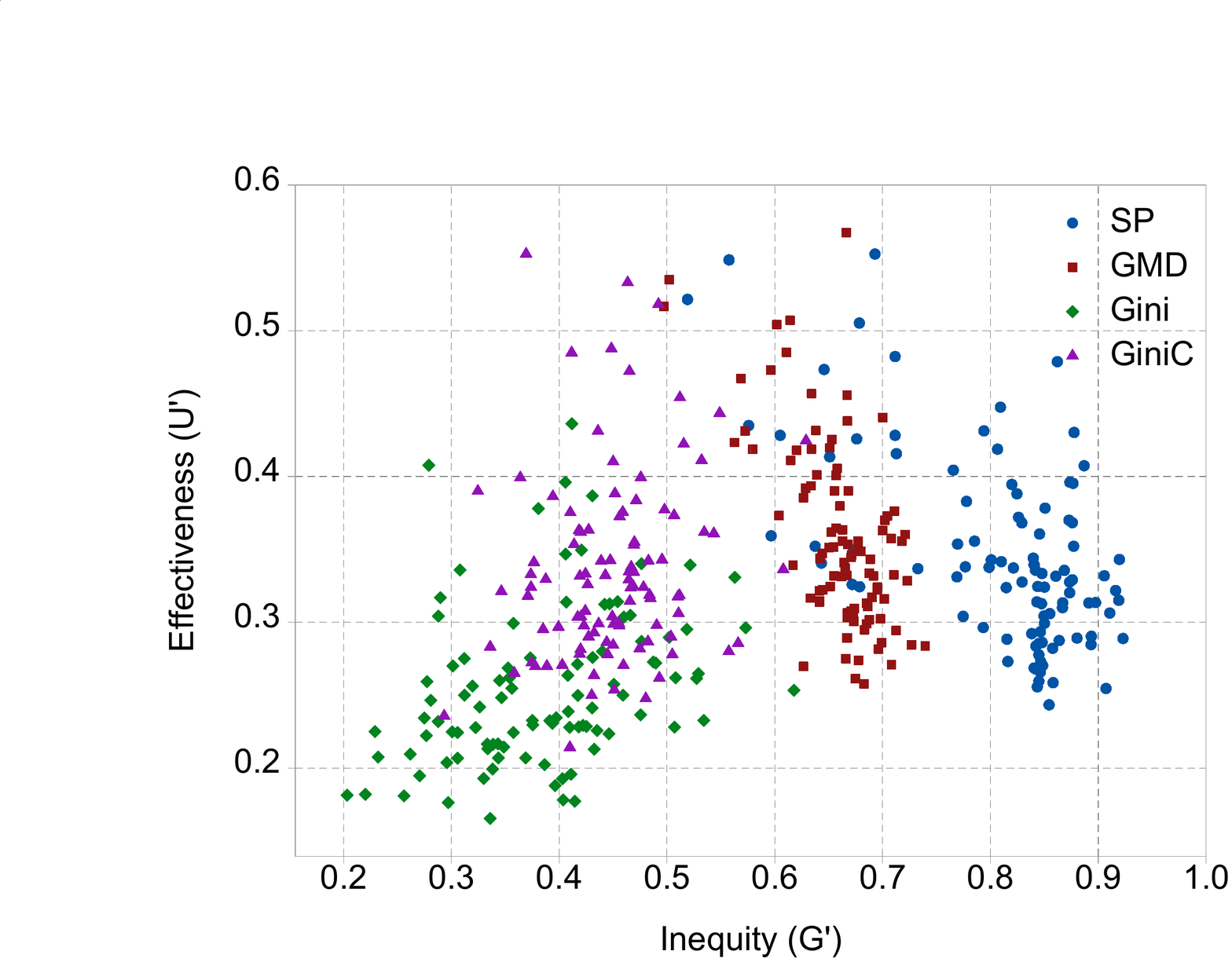}
\vspace{-1.5cm}
\caption{Effectiveness ($U^*_{s'}$) versus inequity ($G^*_{s'}$) for the 100 random realizations.}
  \label{tradeoff1}
\end{figure}

\begin{figure}[H]
\includegraphics[scale=0.6]{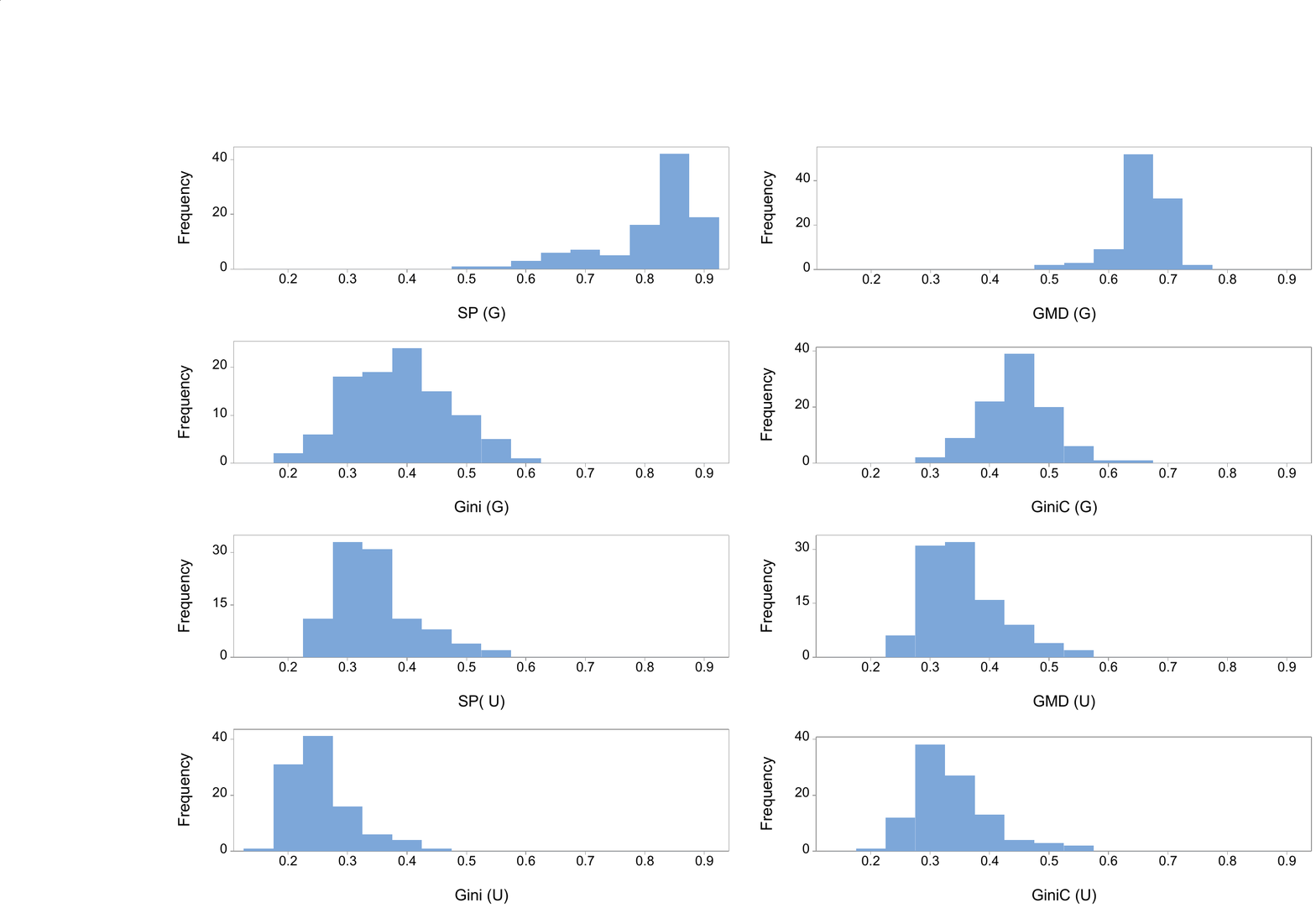}
\vspace{-2.5cm}
\caption{Histograms of $G^*_{s'}$ and $U^*_{s'}$ based on the 100 randomly generated scenarios.}
  \label{Gini_Bound}
\end{figure}

To better visualize the performance of the approaches, the simulation results were summarized into figures and tables.  
Figure \ref{tradeoff1} is a scatter plot of the Gini coefficient ($G^*_{s'}$) versus the effectiveness ($U^*_{s'}$). We can interpret this scatter plot as a Pareto frontier associated with the effectiveness–equity trade-off in which ideal solutions would be in the upper-left quadrant (minimum inequity and maximum effectiveness). Figure \ref{Gini_Bound} shows the empirical distributions of the Gini coefficient and the effectiveness. Finally, Table \ref{tab:payoff} shows the relative benefit of using each approach in terms of inequity and  effectiveness. 

As expected, the SP solutions are concentrated in the upper-right quadrant, which indicates great effectiveness but poor equity. The GMD approach improves equity by 18.58\% on average and maintains considerable effectiveness. Gini and GiniC are clearly dominant in terms of equity. Remarkably, GiniC improves GMD equity by an average of 41.65\% at the expense of a 6.186\% effectiveness reduction. The Gini coefficient histograms indicate that SP (resp. GMD) is narrower around of the mean of 0.8 (resp. 0.66) and has a wider range, relative to Gini and GiniC, achieving values in the range [0.52, 0.92] (resp. [0.50, 0.74]). Left skewness is also clearly visible in SP (resp. GMD), where 96\% (resp. 93\%) of its Gini coefficients are greater than 0.6, which confirms the poor performance in terms of equity. Gini has a wider spread around the mean of 0.39, and achieves the remarkable range of Gini coefficient variation of [0.2,0.62]. Moreover, almost 90\% of the results are less than 0.5, and only 1\% is greater than 0.6. GiniC has a more symmetric histogram and achieves equity levels only slightly worse than Gini (between 0.29 and 0.63). Both Gini and GiniC lead to significantly more near-zero Gini coefficient solutions than the baseline models SP and GMD, which would encourage their adoptions in humanitarian relief chains. The histograms related to effectiveness all show quite similar and expected results. Of course, there is an observable deterioration in effectiveness to enforce greater equity for Gini or GiniC, but GiniC manages to maintain very competitive effectiveness values. This is a telling observation about GiniC, and it can be explained by the fact that it allows a suboptimal portrayal of the Lorenz curve in order to limit the compromise on the effectiveness.

In terms of solvability, it is worth mentioning that all models were well-solved. The elapsed times and optimality gaps for models SP, GMD, Gini, and GiniC were, respectively, 3.63s and 0.00\%, 15.6s and 0.00\%, 3668s and 0.13\%, and 41.1s and 0.06\%. The overall simulation has run for 11.3s, 31.8s, 42,615s, and 35,582s, respectively, and Gini and GiniC had an average optimality gap of (resp.) 2.1\% and 0.02\%. The Gini approach was indeed expected to be harder to solve, mainly because of the additional $|A_s|^2\times |S|$ binary variables required in its formulation. This model can quickly escalate in size for instances with a larger set of disaster-prone areas and scenarios, thus becoming prohibitively expensive without decomposition methods. It is also worth mentioning that the valid inequality \eqref{valI} was fundamental to improve the upper-bound relaxation of Gini (and thus its resolution). The linear relaxation without the valid inequality has an objective function value of 3.4, whereas that with the valid inequality has an objective value of 0.558. Although GiniC has $k_s^2\times|S|$ additional binary variables, it is, as expected, solved much quicker than Gini because $k_s<|A_s|$ for most $s$.  

\vspace{-0.5cm}

\begin{table}[htbp]
\footnotesize
  \centering
  \caption{Relative benefits (\%) in terms of average inequity and effectiveness improvements.}
    
    \begin{tabular}{lrrrrlrrrr}
    \toprule
    Inequity  & SP    & GMD   & Gini  & GiniC & Effectiveness & SP    & GMD   & Gini  & GiniC \\
    \midrule
    SP    & 0     & $-$18.58 & $-$52.49 & $-$44.62 & SP    & 0     & 3.759 & $-$26.26 & $-$2.659 \\
    GMD   &       & 0     & $-$41.65 & $-$31.99 & GMD   &       & 0     & $-$28.93 & $-$6.186 \\
    Gini  &       &       & 0     & 16.56 & Gini  &       &       & 0     & 32.01 \\
    GiniC &       &       &       & 0     & GiniC &       &       &       & 0 \\
    \bottomrule
    \end{tabular}%

    \begin{tablenotes}
		\item \textbf{Note}. Let $\delta[i,]$ be the average metric (inequity or effectiveness) for the model in row $i$ and $\delta[,j]$ be the average metric for the model in column $j$ of cell $(i,j)$. The relative benefit in cell $(i,j)$ in the table is calculated as $\big((\delta[,j]-\delta[i,])/\delta_[i,]\big)\times 100.$
		%The relative benefit is evaluated as follows: $(\delta_i-\delta_j)/\delta_i$, where $\delta_i$ refers to the approach in row $i$ and $\delta_j$ is the approach in column $j$.  
	\end{tablenotes}
  \label{tab:payoff}%
\end{table}%

\normalsize

\section{Conclusion} \label{sec:conclusions}

In this paper, we derive a Gini-coefficient-based objective function for humanitarian logistics planning via the Lorenz curve by means of mixed-integer programming. The Gini coefficient was, in its original form, defined using the Lorenz curve, and our main claim is that existing decision-making approaches in humanitarian logistics use proxies for the true Gini coefficient and thus compromise the equity of the resulting decisions. Through a humanitarian location-allocation problem, we develop two mathematical programming formulations that optimize the effectiveness-equity trade-off. One formulation portrays the Gini coefficient via the optimal (decision-driven) Lorenz curve. The other formulation uses a priori $k$-means clustering to build a cluster-based Lorenz curve, in order to improve the numerical efficiency of the resulting optimization problem. The empirical results based on a realistic case-study of floods and landslides in Rio de Janeiro state, Brazil, reveal that our approaches are promising in reducing inequity in decision-making in comparison to the Gini mean difference, which is one of the most popular proxies for the Gini coefficient employed in humanitarian logistics. Because ``there is no free lunch'', the overall reduction in inequity does come at the cost of lower effectiveness in decision-making, which is acceptable to a certain degree if equitability is the key to having better allocation of resources. We believe that our methodology is general enough to be extended to other optimization problems. In this sense, an interesting direction for future research would be to apply it in other contexts and analyze its performance in comparison to popular equity-based formulations.

\normalsize{
  \bibliography{main}
  \bibliographystyle{spbasic}
}

%\newpage
\appendix 

\section{Input data} \label{ap:input_data}

\begin{landscape}
\begin{table}[htbp]
  \centering
	\scriptsize
	\setlength{\tabcolsep}{3pt}
	\renewcommand{\arraystretch}{1}  
	\caption{Number of victims for different municipalities of the Rio de Janeiro state from 2000 to 2018 \citep{s2id2020}.} \label{tab:sto_victims}%
    \begin{tabular}{crrrrrrrrrrrrrrrrrr}
    \toprule
      City & 2000 & 2001 & 2002 & 2003 & 2004 & 2005 & 2006 & 2007 & 2008 & 2009 & 2010 & 2011 & 2012 & 2013 & 2015 & 2016 & 2017 & 2018 \\
    \midrule
    Teresópolis  & 0   & 10028 & 10763  & 0   &   0 & 0   & 1777 & 2500  & 0   & 112 & 53400 & 49000 & 10162 & 111  & 0   & 102372  & 0   &0  \\
    Petrópolis  & 0   & 10642  & 0   & 1512 & 103372 & 130000  & 0   & 30125 & 45000 & 22200 & 78500 & 70000 & 20000 & 152277 & 104 & 135583  & 0   & 15695 \\
    Nova Friburgo  & 0   &  0  & 0   &  0 & 8 & 2006  & 0   & 80000  & 0   &  0  & 0   & 180000  & 0   &  0  & 0   &  0  & 0   & 0 \\
    São José do Vale do Rio Preto  & 0   &  0  & 0   &  0  & 0   &  0  & 0   & 747  & 0   & 0   & 0   & 20682  & 0   &  0  & 0   &  0  & 0   & 8263 \\
    Bom Jardim  & 0   & 0  & 0   &  0  & 0   &  0  & 0   & 1456  & 0   &   0 & 0   & 12380  & 0   &  0  & 0   &  0  & 0   & 0 \\
    Sumidouro  & 0   &  0  & 0   &  0  & 0   &  0  & 0   & 11000  & 0   &   & 17034 & 20000 & 370  & 0   &  0  & 0   & 0  & 3000 \\
    Areal  & 0   &  0  & 0   &  0  & 0   & 0   & 0   & 206  & 0   &  0 & 130 & 7000  & 0   &   0 & 0   & 0   & 0   & 185 \\
    Santa Maria Madalena  & 0   &  0  & 0   & 0   & 0   &  0  & 0   & 3243 & 6000 & 1204 & 288 & 14049 & 18321  & 0   & 0   & 0   & 185 & 3 \\
    Sapucaia  & 0   &  0  & 0   &  0  & 0   &   0 & 0   & 613  & 0   &   & 5210 & 1520 & 4500  & 0   & 0  & 815  & 0   &  0\\
    São Sebastião do Alto  & 0   &   0 & 0   &  0  & 0   &   0 & 0   & 570  & 0   &  0  & 0   & 8906  & 0   &  0  & 0   &  0  & 0   & 0 \\
    Cordeiro  & 0   &  0  & 0   & 0   & 0   &  0  & 0   & 2724  & 0   &   0 & 0   &   & 1200  & 0   &  0  & 0   &  0 & 0 \\
    Macuco  & 0   &   0 & 0   & 36 & 1000 & 400  & 0   & 784 & 409 & 213  & 0   & 115  & 0   &  0  & 0   &  0  & 0   &  0\\
    Três Rios & 2116  & 0   & 0   & 0   & 57  & 0   & 0  & 1054 & 300 & 35240  & 0   & 2000 & 25000  & 0   &  0  & 0   & 5006 & 0 \\
    \bottomrule
    \end{tabular}%

\vspace{1.5cm}

  \centering
	\scriptsize
	\setlength{\tabcolsep}{3pt}
	\renewcommand{\arraystretch}{1}  
	\caption{Capacities ($\kappa^{\mbox{\scriptsize resp}}_{\ell n}$) and associated fixed cost ($c^{\mbox{\scriptsize o}}_{\ell n}$) for the different size and response facilities locations.} \label{tab:facilities}%
    \begin{tabular}{ccrrrrrrrrrrrrr}
    \toprule
&   &   &   &   &  \multicolumn{1}{c}{São José} &   &   &   & \multicolumn{1}{c}{Santa} &   &   &   &   &  \\
&   &   &   & \multicolumn{1}{c}{Nova} & \multicolumn{1}{c}{do Vale do} & \multicolumn{1}{c}{Bom} &   &   & \multicolumn{1}{c}{Maria} &   & \multicolumn{1}{c}{São Sebastião} &   &   & \multicolumn{1}{c}{Três} \\
& \multicolumn{1}{c}{Level} & \multicolumn{1}{c}{Teresópolis} & \multicolumn{1}{c}{Petropolis} & \multicolumn{1}{c}{Friburgo} & \multicolumn{1}{c}{Rio Preto} & \multicolumn{1}{c}{Jardim} & \multicolumn{1}{c}{Sumidouro} & \multicolumn{1}{c}{Areal} & \multicolumn{1}{c}{Madalena} & \multicolumn{1}{c}{Sapucaia} & \multicolumn{1}{c}{do Alto} & \multicolumn{1}{c}{Cordeiro} & \multicolumn{1}{c}{Macuco} & \multicolumn{1}{c}{Rios} \\
    \midrule
\multirow{4}[2]{*}{$\kappa^{\mbox{\scriptsize resp}}_{\ell n}$} 
& Very large & 114096 & 114096 & 114096 & 34514 & 34514 & 34514 & 34514 & 34514 & 34514 & 34514 & 34514 & 34514 & 34514 \\
& Large & 74123 & 74123 & 74123 & 22434 & 22434 & 22434 & 22434 & 22434 & 22434 & 22434 & 22434 & 22434 & 22434 \\
& Medium & 48180 & 48180 & 48180 & 14582 & 14582 & 14582 & 14583 & 14582 & 14584 & 14582 & 14585 & 14582 & 14586 \\
& Small & 26873 & 39973 & 47250 & 5430 & 3250 & 5250 & 1838 & 4810 & 1368 & 2338 & 716 & 263 & 9251 \\
    \midrule
    \multirow{4}[2]{*}{$c^{\mbox{\scriptsize o}}_{\ell n}$} 
& Very large & 16635197 & 16635197 & 16635197 & 5032141 & 5032141 & 5032141 & 5032141 & 5032141 & 5032141 & 5032141 & 5032141 & 5032141 & 5032141 \\
& Large & 10807133 & 10807133 & 10807133 & 3270877 & 3270877 & 3270877 & 3270877 & 3270877 & 3270877 & 3270877 & 3270877 & 3270877 & 3270877 \\
& Medium & 7024644 & 7024644 & 7024644 & 2126056 & 2126056 & 2126056 & 2126201 & 2126056 & 2126347 & 2126056 & 2126493 & 2126056 & 2126639 \\
& Small & 3918083 & 5828063 & 6889050 & 791694 & 473850 & 765450 & 267980 & 701298 & 199454 & 340880 & 104393 & 38345 & 1348796 \\
    \bottomrule
    \end{tabular}%
\end{table}%
\end{landscape}

\begin{table}[htbp]
	\centering
	\caption{Summary of the relief aid characteristics.}
	\begin{tabular}{lccccc}
		\toprule
		Relief aid     & \mbox{Length in}$^a$ & Coverage$^b$      & Volume in $m^3$  & Prep. capacity$^c$  & Prep. cost \\
		($r$)				& days (\mbox{Length}$_r$)& (\# people) & ($\rho_r$)         & in units ($\theta^{\mbox{\scriptsize max}}_{r}$)  & in BRL ($c^{\mbox{\scriptsize p}}_{rn}$) \\\midrule
		Water           & 7 & 1       & 0.005 & 19,691,500 & 16 \\
		Food            & 1 & 4       & 0.04  & 188,140   & 261 \\
		Hygiene kits    & 1 & 4       & 0.04  & 188,140   & 190 \\
		Cleaning kits   & 1 & 4       & 0.03  & 188,140   & 142 \\
		Dormitory kits  & 1 & 1       & 0.03  & 752,560   & 143 \\
		Mattress        & 1 & 1      & 0.017  & 752,560   & 238 \\ \bottomrule
	\end{tabular}
	\begin{tablenotes}
		\item $^a$ `Length' represents how many days a victims need to be supplied with relief aid $r$ in a horizon of one week. $^b$ `Coverage' shows how many people are covered by one unit of relief aid $r$. $^c$ `Prep. capacity' refers to the maximum quantity of each aid $r$ that could be acquire.
	\end{tablenotes}
	\label{sec:aid}
\end{table}

\begin{landscape}
\begin{table}[htbp]
  \centering
	\scriptsize
	\setlength{\tabcolsep}{3pt}
	\renewcommand{\arraystretch}{1}  
	\caption{Distance (in km) between nodes of the network.} \label{table:distance}%
    \begin{tabular}{lrrrrrrrrrrrrr}
    \toprule
      &   &   &   & \multicolumn{1}{c}{São José} &   &   &   & \multicolumn{1}{c}{Santa} &   &   &   &   &  \\
      &   &   & \multicolumn{1}{c}{Nova} & \multicolumn{1}{c}{do Vale do} & \multicolumn{1}{c}{Bom} &   &   & \multicolumn{1}{c}{Maria} &   & \multicolumn{1}{c}{São Sebastião} &   &   & \multicolumn{1}{c}{Três} \\
      & \multicolumn{1}{c}{Teresópolis} & \multicolumn{1}{c}{Petropolis} & \multicolumn{1}{c}{Friburgo} & \multicolumn{1}{c}{Rio Preto} & \multicolumn{1}{c}{Jardim} & \multicolumn{1}{c}{Sumidouro} & \multicolumn{1}{c}{Areal} & \multicolumn{1}{c}{Madalena} & \multicolumn{1}{c}{Sapucaia} & \multicolumn{1}{c}{do Alto} & \multicolumn{1}{c}{Cordeiro} & \multicolumn{1}{c}{Macuco} & \multicolumn{1}{c}{Rios} \\
    \midrule
    Teresópolis & 0 & 67.6 & 77.1 & 39.4 & 92.4 & 63.2 & 55.4 & 162 & 85 & 150 & 117 & 129 & 80.4 \\
    Petropolis & 67.6 & 0 & 127 & 69.9 & 143 & 113 & 40.1 & 220 & 93.9 & 208 & 168 & 180 & 69.4 \\
    NovaFriburgo & 77.1 & 127 & 0 & 92.8 & 22.5 & 45.4 & 119 & 91.7 & 84.4 & 80.5 & 47.3 & 59.4 & 126 \\
    São Josédo Vale doRio Preto & 39.4 & 69.9 & 92.8 & 0 & 103 & 55.2 & 27.4 & 156 & 34.1 & 144 & 111 & 123 & 47.2 \\
    BomJardim & 92.4 & 143 & 22.5 & 103 & 0 & 51.8 & 129 & 69.5 & 92.2 & 58.3 & 25 & 37.1 & 134 \\
    Sumidouro & 63.2 & 113 & 45.4 & 55.2 & 51.8 & 0 & 81.3 & 105 & 41.2 & 93.5 & 60.3 & 72.4 & 82.8 \\
    Areal & 55.4 & 40.1 & 119 & 27.4 & 129 & 81.3 & 0 & 175 & 49.6 & 164 & 135 & 143 & 24.5 \\
    SantaMariaMadalena & 162 & 220 & 91.7 & 156 & 69.5 & 105 & 175 & 0 & 126 & 22.3 & 49.6 & 34.4 & 167 \\
    Sapucaia & 85 & 93.9 & 84.4 & 34.1 & 92.2 & 41.2 & 49.6 & 126 & 0 & 115 & 85.5 & 93.5 & 50.6 \\
    São Sebastiãodo Alto & 150 & 208 & 80.5 & 144 & 58.3 & 93.5 & 164 & 22.3 & 115 & 0 & 36.4 & 21.2 & 156 \\
    Cordeiro & 117 & 168 & 47.3 & 111 & 25 & 60.3 & 135 & 49.6 & 85.5 & 36.4 & 0 & 17.6 & 127 \\
    Macuco & 129 & 180 & 59.4 & 123 & 37.1 & 72.4 & 143 & 34.4 & 93.5 & 21.2 & 17.6 & 0 & 135 \\
    TrêsRios & 80.4 & 69.4 & 126 & 47.2 & 134 & 82.8 & 24.5 & 167 & 50.6 & 156 & 127 & 135 & 0 \\
    \bottomrule
    \end{tabular}%
\end{table}%
\end{landscape}

\end{document}